\documentclass[12pt, a4paper]{article}

\usepackage[applemac]{inputenc}
\usepackage{amsmath,amssymb,amsfonts,amsthm}
\usepackage{bbm}
\usepackage{enumerate}
\usepackage{graphicx}
\usepackage{geometry}
\usepackage{url}

\newcommand{\Pro}{\mathbf{P}}
\newcommand{\Es}{\mathbf{E}}

\newcommand{\m}[1]{\mathbbm{#1}}

\usepackage[english]{babel}

\theoremstyle{definition}
\newtheorem{thm}{{\it Theorem}}[section]
\newtheorem{Lemma}{{\it Lemma}}[section]

\newtheorem{rem}{{\it Remark}}[section]

\title{Asymptotic properties of the MLE for the autoregressive process coefficients under stationary Gaussian noise}
\author{Alexandre Brouste, Chunhao Cai  and Marina Kleptsyna}
\date{}
\begin{document}

\maketitle
\begin{abstract}
In this paper we are interested in the Maximum Likelihood Estimator (MLE) of the vector  parameter of an autoregressive process of order~$p$ with regular stationary Gaussian noise. We exhibit
the large sample asymptotical properties of the MLE under  very mild conditions. Simulations are done for fractional Gaussian noise (fGn), autoregressive noise (AR(1)) and moving average noise (MA(1)).
\end{abstract}
\section{Statement of the problem}
\subsection{Introduction}
The problem of parametric estimation in classical autoregressive (AR) models generated by  white noises noises has been studied for decades. In particular, for such autoregressive  models of order $1$ (AR(1)) consistency and
many other asymptotic properties (distribution, bias, quadratic error) of the  Maximum Likelihood Estimator (MLE) have been completely analyzed in all possible cases: stable, unstable and explosive (see, \textit{e.g.}, \cite{Anderson1959, ChanWei, Rao1978, Rubin, White1958, White1959}). Concerning autoregressive models of  order~$p$ (AR(p)) with  white noises, the results about the asymptotic behavior of the MLE are less exhaustive but there are still many contributions in the  literature (see, \textit{e.g.},\cite{Anderson1959, duflo,  LaiWei, Alain1, MannWald, Rao1959}).

In the past thirty years numerous papers have been devoted to the statistical analysis of AR processes which may represent long memory phenomenons as encountered in various fields as econometrics \cite{greene1}, hydrology \cite{lawrance}  or biology \cite{maxim}. Of course the relevant models exit from the white noise frame evoked above
 and they involve more or less specific structures of dependence in  the perturbations (see, \textit{e.g.}, \cite{andel, Fox,  granger, hosking, sethuraman, Yajima} for contributions and other references).

General conditions under which the MLE is consistent and asymptotically normal for stationary sequences
have been given in \cite{sweeting}. In order to apply this result, it would be necessary to study the second derivatives of the covariance matrix of the observation sample $(X_1, \ldots, X_N)$. To avoid this difficulty,
some authors followed an other approach suggested by Whittle \cite{Fox} (which is not MLE) for stationary sequences.  But even in autoregressive  models of order $1$  as soon as $|\vartheta| >1$, the process is not stationary anymore and it is not possible to apply theorems in \cite{Fox} to deduce estimator properties.

In the present paper we deal with an AR(p) model  generated by an arbitrary regular stationary Gaussian noise. We exhibit an explicit formula for the MLE of the parameter and we analyze its asymptotic properties.

\subsection{Statement of the problem}

We  consider an  AR(p) process $(X_n,\,n \geq 1)$ defined by  the recursion
\begin{equation}\label{ARp}
X_n=\sum_{i=1}^p\vartheta_iX_{n-i}+ \xi_n, \quad n\geq 1, \quad X_{r}=0,\quad r=0,\,-1, \ldots,\,-(p-1),
\end{equation}
where $\xi=(\xi_n,\,n\in \mathbb{Z})$ is a centered regular stationary Gaussian sequence, {\it i.e.}
\begin{equation}\label{cond spec den}
\int_{-\pi}^{\pi}\left|\ln f_{\xi}(\lambda)\right|d\lambda <\infty,
\end{equation}
where $f_{\xi}(\lambda)$ is the spectral density of $\xi$.
We suppose that the covariance  $c=(c(m,n),\, m,n \ge1)$,
 where
\begin{equation}\label{covar}
    \mathbf{E}\xi_m\xi_n=c(m,\,n)=\rho(|n-m|), \quad \rho(0)=1,
\end{equation}
is positive defined.

For a fixed value of the parameter $\vartheta=(\vartheta_1,\,\dots,\,\vartheta_p) \in\mathbb{R}^p $, let $\mathbf{P}_{\vartheta}^N$ denote the probability measure induced by $X^{(N)}$.  Let $\mathcal{L}(\vartheta,\,X^{(N)})$ be the likelihood function defined by the Radon-Nikodym derivative of $\mathbf{P}_{\vartheta}^N$ with respect to the Lebesgue measure.
Our goal is to study the large sample asymptotical properties of the Maximum Likelihood Estimator (MLE) $\widehat{\vartheta}_N$ of $\vartheta$ based on the observation sample $X^{(N)}=\left(X_1, \ldots, X_N\right)$:
\begin{equation}\label{def MLE}
\hat{\vartheta}_N=\sup_{\vartheta \in  \mathbb{R}^p} \mathcal{L}(\vartheta,\,X^{(N)}).
\end{equation}

At first, preparing for the analysis of the consistency (or strong consistency) of $\hat{\vartheta}_N$ and its  limit distribution we transform our observation model into an "equivalent" model with independent Gaussian noises. This allows to write explicitly the MLE and
actually, the difference between $ \hat{\vartheta}_N$  and the real value $\vartheta$ appears as  the product of  a  martingale by  the inverse of its bracket process. Then we can use Laplace transforms computations to prove the asymptotical properties of the MLE.

The paper is organized as follows. Section~ \ref{Results} contains  theoretical results and simulations. Sections~ \ref{prelim} and \ref{Auxilres} are devoted to preliminaries and auxiliary results. The proofs of the main results are presented in Section~ \ref{proofs}.

\paragraph{Acknowledgments} We would like to thank Alain Le Breton for
very fruitful discussions and his interest for this work.

\section{ Results and illustrations}\label{Results}
\subsection{Results}
We define the $p\times p$ companion matrix $A_0$ and the  vector $b\in \mathbb{R}^p$ as follows:
\begin{equation}\label{def A_0}
A_0=\left(
      \begin{array}{ccccc}
        \vartheta_1 & \vartheta_2 & \cdots & \vartheta_{p-1} & \vartheta_p \\
        1 & 0 & \cdots & 0 & 0 \\
        0 & 1 & \cdots & 0 & 0 \\
        \vdots & \vdots & \ddots & \vdots & \vdots \\
        0 & 0 & \cdots & 1 & 0 \\
      \end{array}
    \right)
, \quad
b=\left(
    \begin{array}{c}
      1 \\
      \mathbf{0}_{(p-1)\times 1} \\
    \end{array}
  \right).
\end{equation}
Let  $r(\vartheta)$ be the spectral radius of  $A_0$. The  following results hold:
\begin{thm}\label{p-dimension}
Let $p \geq 1$ and  the parameter set be:
 \begin{equation}\label{stable}
\Theta=\left\{\vartheta\in \mathbb{R}^p \, | \,r(\vartheta)<1\right\}.
\end{equation}
The MLE $ \hat{\vartheta}_N$ is consistent, {\it i.e.}, for any $\vartheta \in \Theta$ and $\nu >0$,
\begin{equation}
\lim_{N\rightarrow \infty} \Pro_\vartheta^{N} \left\{\left\| \hat{\vartheta}_N-\vartheta \right\|> \nu \right\}=0\,,
\end{equation}
and  asymptotically normal
\begin{equation}\label{eq:asympnormal}
\sqrt{N}\left(\hat{\vartheta}_N-\vartheta\right)\overset{\textit{law}}{\Rightarrow}\mathcal{N}(\mathbf{0},\,\mathcal{I}^{-1}(\vartheta)),
\end{equation}
where $\mathcal{I}(\vartheta)$ is the unique solution of  the Lyapounov  equation:
\begin{equation}\label{Lyapunov}
\mathcal{I}(\vartheta)=A_0\mathcal{I}(\vartheta)A_0^*+bb^*,
\end{equation}
for $A_0$ and  $b$  defined in \eqref{def A_0}.

Moreover  we have the  convergence of the moments: for any $\vartheta \in \Theta$ and $q>0$
\begin{equation}
\lim_{N\rightarrow \infty}   \left| \Es_\vartheta \left\|\sqrt{N} \left( \hat{\vartheta}_N-\vartheta\right) \right \|^q-  \Es \left\|\eta\right\|^q \right| =0,
\end{equation}
where $\|\,\|$ denotes the  Euclidian norm on $\mathbb{R}^p$ and  $\eta$ is a zero mean Gaussian random vector with covariance matrix ${\cal I}(\vartheta)^{-1}.$
\end{thm}
\begin{rem}\label{comp_clas_case}
It is worth to emphasize that the asymptotic covariance $\mathcal{I}^{-1}(\vartheta)$ is actually the same as in the standard case where $(\xi_{n})$ is a white noise.  (\textit{cf}. \cite{?}).
\end{rem}

In the case $p=1$ we can strengthen the assertions of Theorem~\ref{p-dimension}. In particular, the strong consistency and uniform convergence on compacts of the moments hold.
%

\begin{thm}\label{strong}
Let $p=1$ and the parameter set be $\Theta=\mathbb{R}.$ The MLE $ \hat{\vartheta}_N$ is strongly consistent, \textit{i.e.} for any $\vartheta \in \Theta $
\begin{equation}\label{strong_consist}
\lim_{N\rightarrow \infty}\hat{\vartheta}_N=\vartheta \quad \quad   a.s..
\end{equation}
Moreover, $ \hat{\vartheta}_N$ is uniformly consistent and satisfies the uniform convergence of the moments on compacts $\m{K} \subset (-1,\,1)$, {\it i.e.} for any $\nu >0:$
\begin{equation}
\lim_{N\rightarrow \infty} \sup_{\vartheta \in \mathbb{K}} \Pro_\vartheta^{N} \left\{\left| \hat{\vartheta}_N-\vartheta \right|> \nu \right\}=0\,,
\end{equation}
and  for any $q>0:$,
\begin{equation}
\lim_{N\rightarrow \infty} \sup_{\vartheta \in \mathbb{K}}  \left| \Es_\vartheta \left|\sqrt{N} \left( \hat{\vartheta}_N-\vartheta\right) \right|^q-  \Es \left|\eta\right|^q \right| =0,
\end{equation}
where $\eta \sim \mathcal{N}(0,\,1- \vartheta^2).$
\end{thm}
\begin{rem}\label{spect_dens_cov}
It is worth mentioning that condition \eqref{cond spec den} can be rewritten in terms of the covariance function $\rho$ : $\rho(n)\sim n^{-\alpha},\, \alpha>0$.
\end{rem}

\subsection{Simulations}\label{sec:simu}
In this section we present for $p=1$ three illustrations of the behavior of the MLE  corresponding to noises which are MA(1), AR(1) and fGn.

\paragraph{Moving average noise MA(1) }

Here we consider MA(1) noises where
$$ \xi_{n+1} = \frac{1}{\sqrt{1+\alpha^2}} (\varepsilon_{n+1} + \alpha \varepsilon_{n} ), \quad n\geq 1, $$
where $( \varepsilon_n, n\geq 1)$ is a sequence of i.i.d. zero-mean standard Gaussian variables.
Then the  covariance function is given by
$$
\rho(|n-m|)= \mathbbm{1}_{\{|n-m|=0\}} + \frac{\alpha}{1+\alpha^2} \mathbbm{1}_{\{|n-m|=1\}}.
$$
Condition \eqref{cond spec den} is fulfilled for $|\alpha| <1$.

\paragraph{ Autoregressive noise (AR(1))}
Here  we consider stationary autoregressive AR(1) noises where
$$ \xi_{n+1} = \sqrt{1-\alpha^2}\varepsilon_{n+1} + \alpha \xi_{n}, \quad n\geq 1,$$
where $( \varepsilon_n, n\geq 1)$ is a sequence of i.i.d. zero-mean standard Gaussian variables.
Then the  covariance function is
$$
\rho(|n-m|)= \alpha^{|n-m|}.
$$
Condition \eqref{cond spec den} is fulfilled for $|\alpha| <1$.

\paragraph{ Fractional Gaussian noise fGn}
Here the covariance function of $(\xi_{n})$ is
$$
\rho(|m-n|)= \frac{1}{2}\left(|m-n+1|^{2H} -2|m-n|^{2H} + |m-n-1|^{2H}\right),
$$
for a known Hurst exponent $H\in (0,1).$  For simulation of the fGn we use Wood and Chan method~(see~\cite{fGnsim2}).
The explicit formula for the spectral density of fGn sequence has been exhibited  in \cite{sinai}.
Condition \eqref{cond spec den} is fulfilled for any $H\in (0,1). $

\vskip 12pt

\begin{figure}[ht]
\includegraphics[width=6cm, height=5cm]{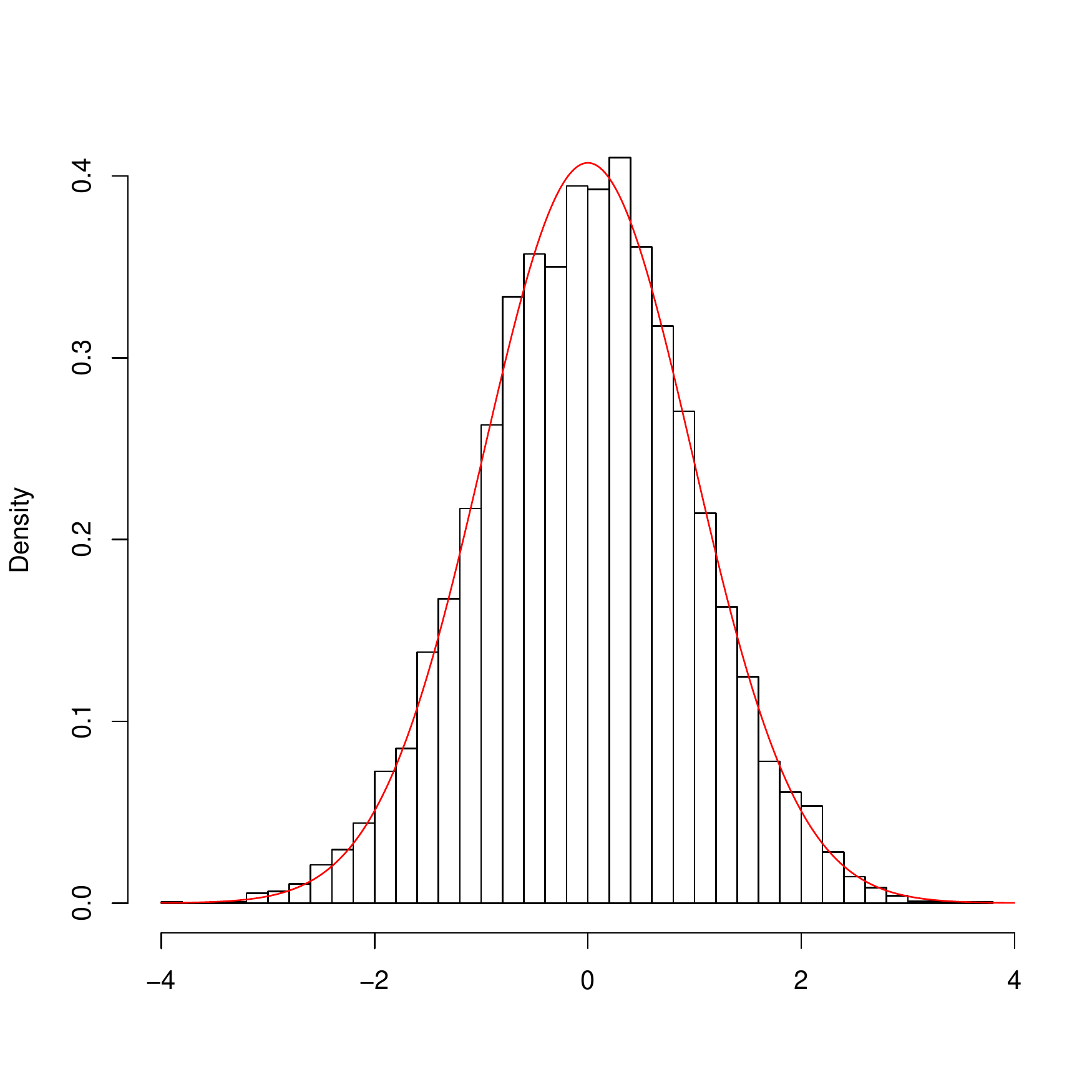}\hfill
\includegraphics[width=6cm, height=5cm]{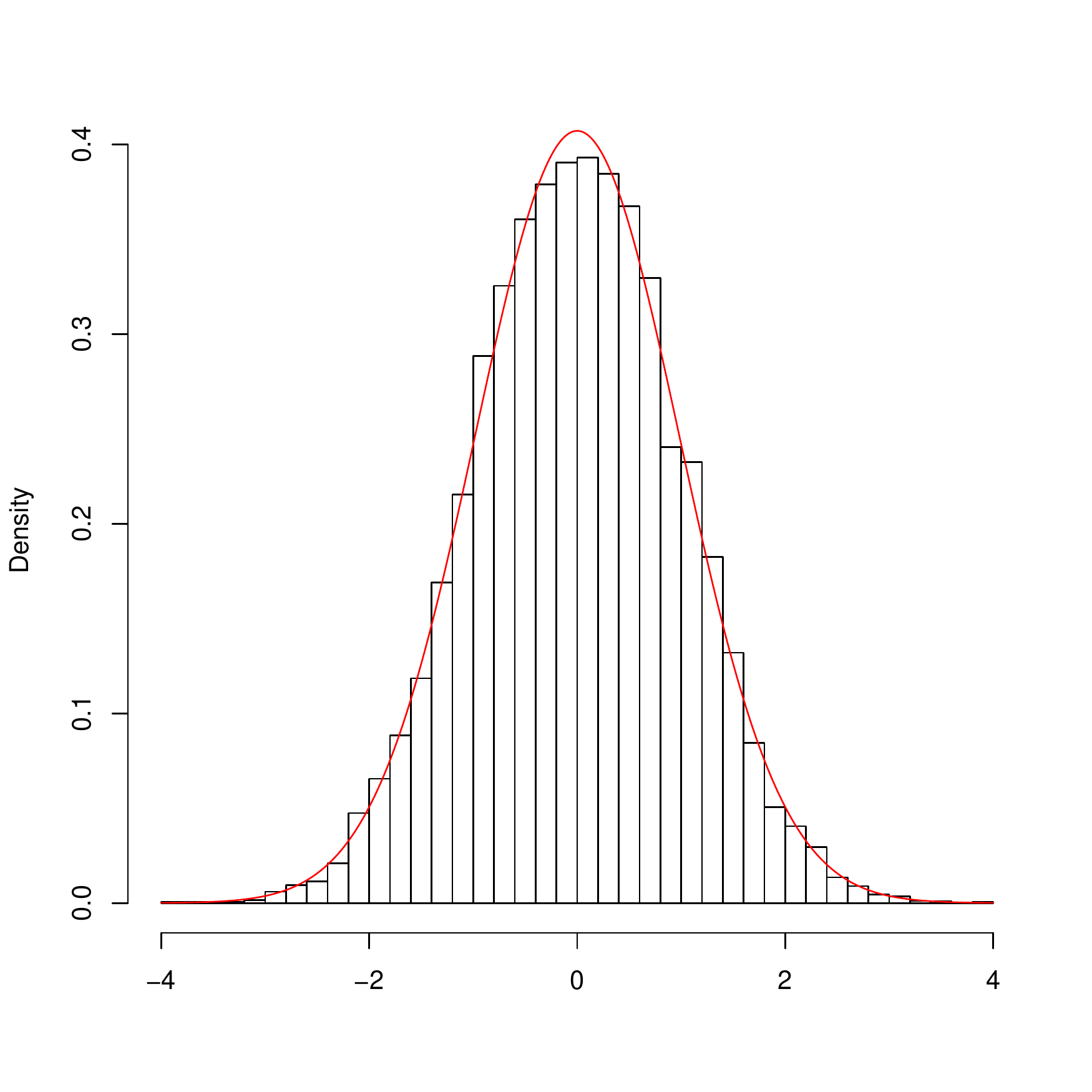}
\includegraphics[width=6cm, height=5cm]{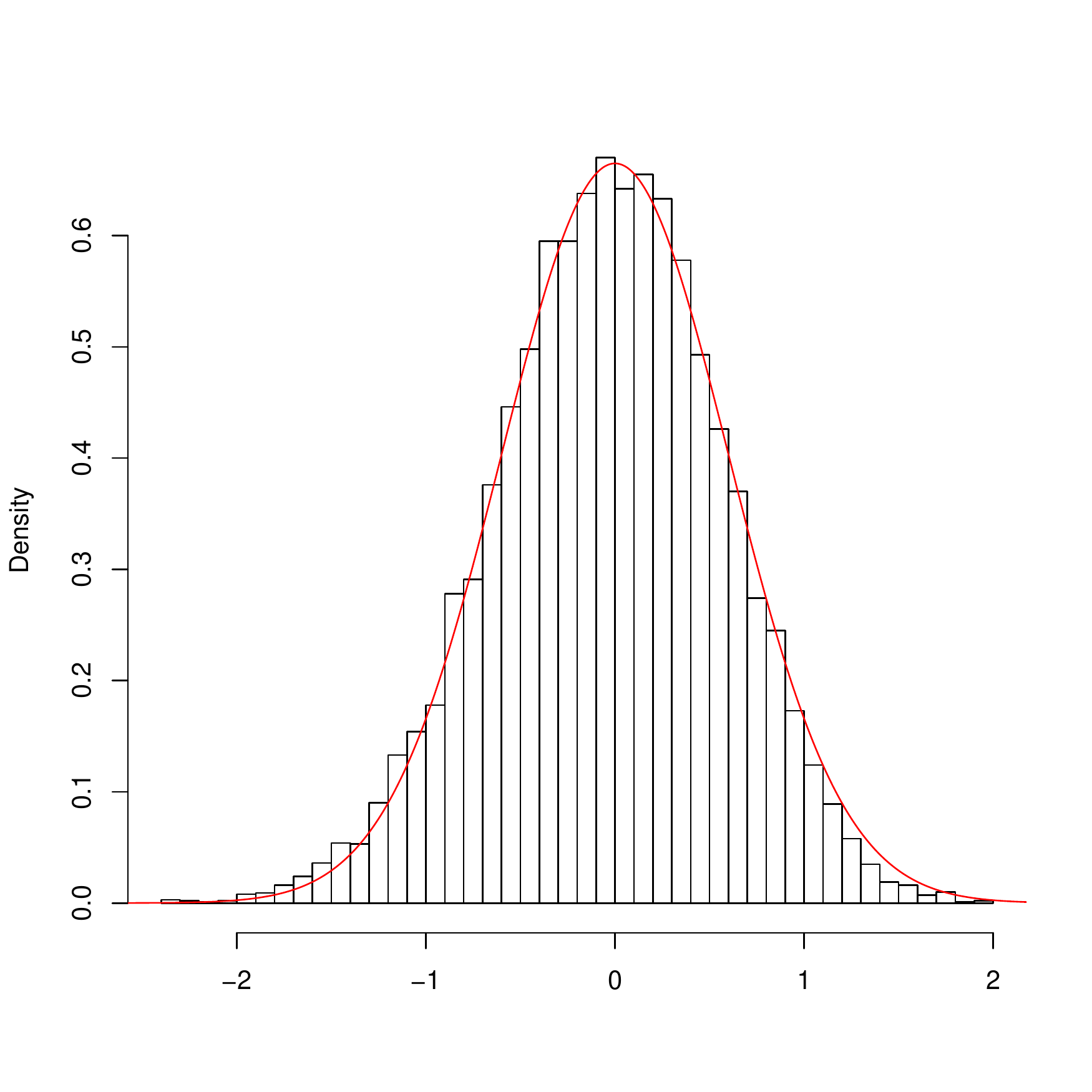}\hfill
\includegraphics[width=6cm, height=5cm]{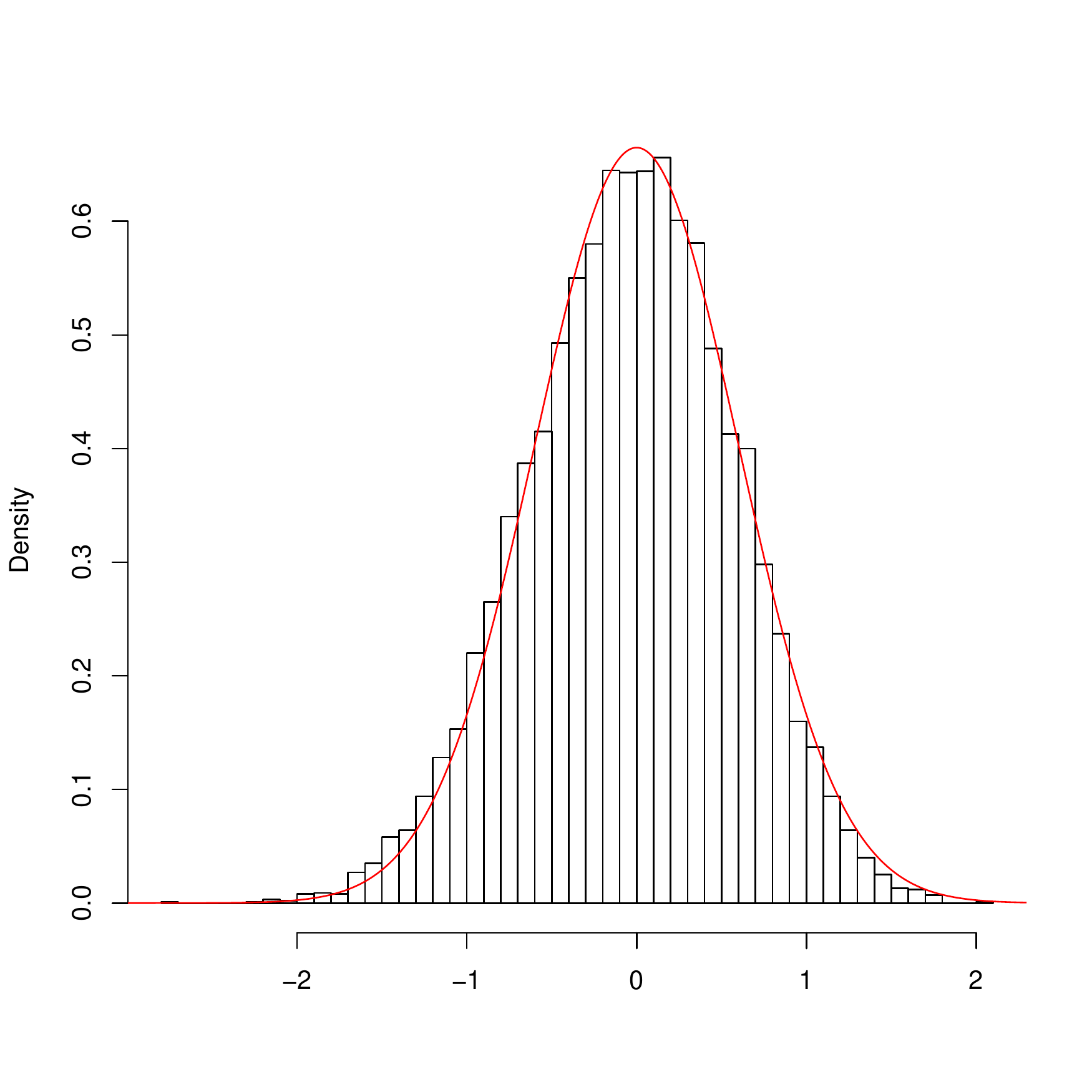}

\caption{Asymptotical normality $N=2000$ for the MLE in different cases by Monte-Carlo simulation of $M=10000$ independent replications for 
AR(1) noises (top left) and  MA noises (top right), both for $\alpha=0.4$ and $\vartheta=0.2$, and fGn noises for $H=0.2$ (bottom left) and for $H=0.8$ (bottom right) both for $\vartheta=0.8$.}\label{fig:sim}
\end{figure}

On Figure  \ref{fig:sim} we can see  that  in conformity with Theorem \ref{strong},   in the three cases the MLE is asymptotically normal with the same limiting variance as in the classical i.i.d. case.

\section{Preliminaries}\label{prelim}
\subsection{Stationary Gaussian sequences}
We begin with some well known properties of a stationary scalar Gaussian sequence $\xi=(\xi_{n})_{n\geq 1}$. We denote
by $(\sigma_n \varepsilon_{n})_{n\geq 1}$ {\em the innovation type sequence } of $\xi$ defined by
$$ \sigma_1\varepsilon_1=\xi_1,  \quad \sigma_n \varepsilon_{n}=\xi_{n}- \Es(\xi_{n}\,|\,\xi_{1}, \dots, \xi_{n-1}), \quad n \geq 2,$$
where  $\varepsilon_{n}\sim
\mathcal{N}(0,1)$, $n\geq 1$ are independent. It follows from the Theorem of Normal Correlation (\cite{Liptser}, Theorem 13.1) that there exists  a  deterministic kernel
denoted by
$k(n,m)$, $n\geq 1$, $m\leq n$, such that

\begin{equation}\label{k}
\sigma_n \varepsilon_{n}=\sum_{m=1}^{n} k(n,m) \xi_{m},\quad k(n,n)=1.
\end{equation}
In the sequel, for $n\geq 1$, we  denote by $\beta_{n-1}$ the partial correlation coefficient
\begin{equation}\label{correlation beta}
-k(n,1)=\beta_{n-1},\, n \geq 1.
\end{equation}
The following relations between $k(\cdot,\,\cdot)$, the covariance function $\rho(\cdot)$ defined by \eqref{covar}, the sequence of partial correlation coefficients  $\left(\beta_n\right)_{n \geq 1}$ and the variances of innovations $\left(\sigma_n^2\right)_{n \geq 1}$  hold (see Levinson-Durbin algorithm \cite{durbin})
\begin{equation}\label{sigma}
\sigma_{n}^{2} = \prod_{m=1}^{n-1} ( 1 - \beta_m^2 ), \quad n\geq 2, \quad \sigma_1=1,
\end{equation}

\begin{equation}\label{correlation}
\sum_{m=1}^nk(n,\,m)\rho(m)=\beta_n\sigma_n^2,
\end{equation}

\begin{equation}\label{func k}
k(n+1,\,n+1-m)=k(n,\,n-m)-\beta_nk(n,\,m).
\end{equation}
Since we assume the positive definiteness of the covariance $c(\cdot,\,\cdot)$, there also exists an inverse deterministic kernel $K=(K(n,\,m),\,n\geq 1,\,m\leq n)$ such that
\begin{equation}\label{repdir}
\xi_n = \sum_{m=1}^{n} K(n,m) \sigma_m \varepsilon_{m}.
\end{equation}

\begin{rem}\label{choldecom}
Actually, kernels $k$ and $K$ are nothing but the ingredients of the Choleski decomposition of  covariance and inverse of covariance matrices. Namely,
$$
\Gamma_{n}^{-1}=k_{n}D_{n}^{-1}k_{n}^{*} \quad \mbox{and} \quad\Gamma_{n}=K_{n}^{*}D_{n}K_{n},
$$
where $\Gamma_{n}=((\rho(|i-j|)))$ , $k_{n}$ and $K_{n}$ are  $n\times n$ lower triangular matrices with ones as diagonal entries and $k(i,j)$ and $K(i,j)$ as subdiagonal entries respectively
and  $D_n$  is an $n\times n$ diagonal matrix with $\sigma_{i}^{2}$ as diagonal entries. Here $*$ denotes  the transposition.
\end{rem}

\begin{rem}\label{betta_l2}
It is worth mentioning that condition \eqref{cond spec den} implies that
\begin{equation}\label{sommab_betta}
\sum_{n\ge1}\beta_{n}^{2}<\infty.
\end{equation}
Indeed,
 for every regular stationary Gaussian sequence $\xi=(\xi_n,\,n\in \mathbb{Z})$, there exists a sequence of i.i.d $\mathcal{N}(0,1)$ random variables $(\widetilde{\varepsilon}_n,\,n\in \mathbb{Z})$ and a sequence of real numbers $a_k,\,k\geq 0$  with  $a_0\neq 0$ such that:
$$
\xi_n=\sum_{k=0}^{\infty}a_k\widetilde{\varepsilon}_{n-k},
$$
 and for all $n\in \mathbb{Z}$ the $\sigma$-algebra generated by  $(\xi_k)_{-\infty <k\le n }$ coincides with the $\sigma$-algebra generated by $(\widetilde{\varepsilon}_k)_{-\infty <k\le n }.$

Note that the variance $\sigma_n^{2}$ of the innovations  is also the one step predicting error and the following
equalities hold thanks to the stationarity of $\xi$:
$$
\lim_{n\rightarrow\infty}\prod_{m=1}^{n-1} ( 1 - \beta_m^2 )=\lim_{n\rightarrow\infty}\sigma_n^{2}
$$
$$
=\lim_{n\rightarrow\infty}\mathbf{E}\left(\xi_n- \mathbf{E}(\xi_n |\xi_1, \cdots \xi_{n-1})\right)^2=
\lim_{n\rightarrow\infty}\mathbf{E}\left(\xi_0- \mathbf{E}(\xi_0|\xi_{-1}, \cdots \xi_{-n+1})\right)^2
$$
$$
=\mathbf{E}\left(\xi_0- \mathbf{E}(\xi_0|\xi_s, s\leq -1)\right)^2=\mathbf{E}\left(\xi_0- \mathbf{E}(\xi_0|\varepsilon_s, s\leq -1)\right)^2= a_0^2>0
$$
which implies \eqref{sommab_betta}.

\end{rem}

\subsection{Model Transformation}\label{model trans}
As usual, for the first step we extend the dimension of the observations in order to work with a first order autoregression in $\mathbb{R}^p$. Namely, let  $Y_n,\, n\ge 1,$ be
$Y_n=( X_n,\, X_{n-1},\,\dots,\,X_{n-(p-1)} )^{*}$
then $Y=(Y_n,\,n\geq 1)$ satisfies the first order autoregressive equation:
\begin{equation}\label{Yn}
Y_n=A_0 Y_{n-1}+b \xi_n,\quad n \geq 1, \quad   Y_0=\mathbf{0}_{p\times 1},
\end{equation}
where $A_0$ and $b$ are defined in \eqref{def A_0}. For the second step we take an appropriate  linear transformation of $Y$ in order to have  i.i.d. noises in the corresponding observations. For this goal
let us introduce the process $Z=\left(Z_n, n \geq 1\right)$ such that
\begin{equation}\label{def Z}
Z_n=\sum_{m=1}^nk(n,\,m)Y_m,   \quad n \geq 1,
\end{equation}
where $k=(k(n,\,m),\,n\geq 1,\,m\leq n)$ is the kernel appearing in \eqref{k}.
Since we have also
\begin{equation}\label{YZ}
Y_n = \sum_{m=1}^{n} K(n,m)  Z_{m},
\end{equation}
where $K=(K(n,\,m),\,n\geq 1,\,m\leq n)$ is the inverse kernel of $k$ (see \eqref{repdir}), the filtration of $Z$ coincides with the filtration of $Y$ (and also the filtration of $X$).
Actually, it was shown in \cite{BK11} that $Z$ can be considered as the first component of a $2p$ dimensional $AR(1)$ process $ \zeta=( \zeta_n,\, n\ge 1)$ governed by  i.i.d. noises. More precisely,
the  process $ \zeta=( \zeta_n,\, n\ge 1)$  defined by :
$$
 \zeta_n=\begin{pmatrix}
Z_n\\
 \displaystyle\sum_{r=1}^{n-1}\beta_r Z_r
\end{pmatrix},
$$
  is a 2p-dimensional Markovian process which satisfies the following equation:
\begin{equation}\label{zeta}
\zeta_{n}=\mathbf{A}_{n-1}\zeta_{n-1}+\ell \sigma_{n}\varepsilon_{n},\quad \,n \geq 1,\, \quad \zeta_0=\mathbf{0}_{2p\times 1},
\end{equation}
where
\begin{equation}\label{def An}
\mathbf{A}_n=\left(
      \begin{array}{cc}
        A_0 & \beta_n A_0 \\
        \beta_n \mathbf{Id}_{p\times p} & \mathbf{Id}_{p\times p} \\
      \end{array}
    \right), \quad \ell=\left(
                          \begin{array}{c}
                            1 \\
                            \mathbf{0}_{(2p-1)\times 1} \\
                          \end{array}
                        \right),
\end{equation}
and $(\varepsilon_n,\, n\geq 1)$ are i.i.d. zero mean standard Gaussian variables.
Now the initial estimation problem is  replaced by the problem of estimation of the unknown parameter $\vartheta$ from the observations of $ \zeta=( \zeta_n,\, n\ge 1)$.

\subsection{Maximum Likelihood Estimator}
It follows directly from equation \eqref{zeta} that  the log-likelihood function is nothing but:
$$
\ln \mathcal{L}(\vartheta,\,X^{(N)})=-\frac{1}{2}\sum_{n=1}^N\left(\frac{\ell^{*}(\zeta_n-\mathbf{A}_{n-1}\zeta_{n-1})}{\sigma_n}\right)^2-\frac{N}{2}\ln 2\pi- \frac{1}{2}\sum_{n=1}^{N}\ln \sigma_n^2
$$
and that the maximum likelihood estimator  $\hat{\vartheta}_N$ is:
\begin{equation}
\hat{\vartheta}_N=\left(\sum_{n=1}^N\frac{a_{n-1}^*\zeta_{n-1}\zeta_{n-1}^*a_{n-1}}{\sigma_n^2}\right)^{-1}\cdot \left(\sum _{n=1}^N\frac{a_{n-1}^*\zeta_{n-1}\ell^{*}\zeta_{n}}{\sigma_n^2}\right).
\end{equation}
Then we can write
\begin{equation}\label{difference}
\hat{\vartheta}_N- \vartheta= 
 \left(\langle M\rangle_N\right)^{-1} \cdot M_N,
\end{equation}
where 
\begin{equation}\label{martbrec}
M_N=\sum _{n=1}^N\frac{a_{n-1}^*\zeta_{n-1}}{\sigma_n}\varepsilon_n\,, \quad  \quad
\langle M\rangle_N=\sum_{n=1}^N\frac{a_{n-1}^*\zeta_{n-1}\zeta_{n-1}^*a_{n-1}}{\sigma_n^2},
\end{equation}
with
\begin{equation}\label{an}
a_n=\left(
             \begin{array}{c}
               \mathbf{Id}_{p\times p} \\
               \beta_n \mathbf{Id}_{p\times p}
             \end{array}
           \right).
\end{equation}
Note that $(M_{n},\, n\ge 1)$ is a martingale and $(\langle M\rangle_n,\, n\ge 1)$ is its bracket process.
\begin{rem}\label{comp_iid}
It is worth mentioning that in the classical i.i.d. case, \textit{i.e.}, when $\beta_{n}=0,\, n\ge 1,$ $M_{N}$ and $\langle M\rangle_N$ in equations \eqref{difference}-\eqref{martbrec} reduce to:
$$
M_N=\sum _{n=1}^N Y_{n-1}\varepsilon_n \,, \quad  \quad
\langle M\rangle_N=\sum_{n=1}^N Y_{n-1} Y_{n-1}^{*}.
$$
Of course, under the condition $r(\vartheta)<1$ due to the law of the large numbers and the central limit theorem for martingales the following convergences hold:
\begin{equation}\label{conviid}
\mathbf{P}_{\vartheta}-\lim_{N\rightarrow \infty}\frac{1}{N}\langle M\rangle_N=\mathcal{I}(\vartheta),\, \frac{1}{\sqrt{N}}M_N \overset{\textit{law}}{\Rightarrow} \mathcal{N}\left(\mathbf{0},\,\mathcal{I}(\vartheta)\right),
\end{equation}
where $\mathcal{I}(\vartheta)$ is the unique solution of  the Lyapounov  equation \eqref{Lyapunov}. This implies immediately the consistency and the asymptotic normality   of the  MLE.
\end{rem}

\section{Auxiliary results}\label{Auxilres}
Actually, the proof of Theorems \ref{p-dimension}- \ref{strong} is crucially based on the asymptotic study for $N$ tending to infinity of the Laplace transform:
\begin{equation}\label{Ln}
L_N^{\vartheta}(\mu)=\mathbf{E}_{\vartheta}\exp\left(-\frac{\mu}{2}\alpha^*\langle M\rangle_{_{N}}\alpha\right),
\end{equation}
for arbitrary $\alpha \in \mathbb{R}^{p}$ and a positive real number $\mu$, where $\langle M\rangle_N$ is defined by \eqref{martbrec}. It  can be rewritten as
\begin{equation}\label{Lapzeta}
L_N^{\vartheta}(\mu)=
\Es_{\vartheta} \exp \left(- \frac{\mu}{2} \sum_{n=1}^{N} \zeta_n^* {\cal M}_{n}\zeta_{n} \right),
\end{equation}
where
$\mathcal{M}_n=\displaystyle{\frac{1}{\sigma_{n+1}^2}}a_n\alpha\alpha^{*}a_n^{*}$, $a_{n}$ is defined by \eqref{an} and $\zeta$ satisfies the equation \eqref{zeta}.
In the sequel we will suppose that
 all the eigenvalues of $A_{0}$ are simple and different from $0$. Actually, it is not a real restriction, since the general case can be studied by using small perturbations arguments.
\begin{Lemma}\label{Lapltrexpl}
 The Laplace transform $L_N^{\vartheta}(\mu)$ can be written explicitly in the  following form:
\begin{equation}\label{Lapl_Psi}
L_N^{\vartheta}(\mu)= \left( \left( \prod_{n=1}^{N}  \det  \mathbf{A}_n \right) \det \Psi^1_N \right)^{-\frac{1}{2}},
\end{equation}
where $ \mathbf{A}_{n}$ is defined by equation \eqref{def An} and
\begin{equation}\label{PsiN}
\sigma_{N+1}^2\Psi_N^1=\Psi_{0} \mathbf{J} \prod_{n=1}^{N}(\mathcal{A}_{\mu}\otimes A_1^n+\mathbf{Id}_{2p\times 2p}\otimes A_2^n) \mathbf{J^{*}}\Psi_{0}^{*}.
\end{equation}
Here $\otimes$ is the Kronecker product, $\Psi_{0}=(\mathbf{Id}_{2p\times 2p} \quad \mathbf{0}_{2p\times 2p})$,
\begin{equation}\label{Amu}
\mathcal{A}_{\mu}=\left(
              \begin{array}{cc}
                A_0^{-1} & A_0^{-1}bb^* \\
               \mu  \alpha\alpha^* & A_0^* +\mu \alpha\alpha^*A_0^{-1}bb^*\\
              \end{array}
            \right)
\end{equation}
and
$2\times2$ matrices  $A_1^n,\,A_2^n$  are defined by

\begin{equation}\label{A1n}
A_1^n=\left(
        \begin{array}{cc}
          1 & 0 \\
          -\beta_n & 0 \\
        \end{array}
      \right), \quad \quad A_2^n=\left(
                                  \begin{array}{cc}
                                    0 & -\beta_n \\
                                    0 & 1 \\
                                  \end{array}
                                \right).
\end{equation}

\end{Lemma}
\begin{proof}
 The following equality  can be proved by using the same arguments as those used in  \cite{AMM} (see equations (15) and (27)):

$$
 L_N^{\vartheta}(\mu) = \left( \left( \prod_{n=1}^{N}  \det  \mathbf{A}_n \right) \det \Psi^1_N \right)^{-\frac{1}{2}},
$$
where $\Psi=((\Psi^1_n, \Psi^2_n),\, n\ge 1)$ is the solution of the following equation:

\begin{equation}\label{fpr8}\nonumber
\left\{
\begin{array}{rrrrrrrl}
\Psi^1_{n-1}&=& \Psi^1_{n} \mathbf{A}_{n} &-& \mu  \Psi_{n-1}^2 {\cal M}_{n-1}   , &\quad  n \geq 1,
\\
\Psi^2_{n}&=& \Psi^1_{n} \ell \ell^* \sigma_{n}^2  &+&\Psi^2_{n-1}  \mathbf{A}_{n}^*  ,& \quad n \geq 1,
\end{array}
\right.
\end{equation}
with the initial condition $\Psi_{0}=(\mathbf{Id}_{2p\times 2p} \quad \mathbf{0}_{2p\times 2p})$.
This equation can be rewritten as

\begin{equation}\label{fpr} \nonumber
\left\{
\begin{array}{rrrrrrrl}
\Psi^1_{n}&=& \Psi^1_{n-1} \mathbf{A}_{n}^{-1} &+& \mu  \Psi_{n-1}^2 {\cal M}_{n-1} \mathbf{A}_{n}^{-1}   , &\quad  n \geq 1,
\\
\Psi^2_{n}&=& \Psi^1_{n-1} \mathbf{A}_{n}^{-1} \ell \ell^* \sigma_{n}^2  &+&\Psi^2_{n-1}   \left( \mu  {\cal M}_{n-1}  \mathbf{A}_{n}^{-1} \ell \ell^* \sigma_{n}^2 +  \mathbf{A}_{n}^*  \right)  ,& \quad n \geq 1.
\end{array}
\right.
\end{equation}
Now let us  denote by
$\widetilde{\Psi}_n^1=\sigma_{n+1}^2\Psi_n^1$ and $\widetilde{\Psi}_n^2=\Psi_n^2\left(
\begin{array}{cc}
 \mathbf{Id}_{p\times p} & \mathbf{0}_{p\times p} \\
  \mathbf{0}_{p\times p}& -\mathbf{Id}_{p\times p} \\
  \end{array}
   \right)
$. Then  $\widetilde{\Psi}_n=(\widetilde{\Psi}_n^1\quad \widetilde{\Psi}_n^2)$ satisfies for $n\geq 1$ the following equation
$$
\widetilde{\Psi}_{n}=\widetilde{\Psi}_{n-1}\left(
 \begin{array}{cccc}
  A_0^{-1} & -\beta_{n}\mathbf{Id}_{p\times p} & A_0^{-1}bb^* & \mathbf{0}_{p\times p} \\
   -\beta_{n} A_0^{-1} & \mathbf{Id}_{p\times p} & -\beta_{n}A_0^{-1}bb^* & \mathbf{0}_{p\times p} \\
    \mu\alpha\alpha^*A_0^{-1} & \mathbf{0}_{p\times p} & \mu\alpha\alpha^*A_0^{-1}bb^*+A_0^* & -\beta_{n}\mathbf{Id}_{p\times p} \\
     -\beta_{n}(\mu\alpha\alpha^*A_0^{-1}) & \mathbf{0}_{p\times p} & -\beta_{n}(\mu\alpha\alpha^*A_0^{-1}bb^*+A_0^*) & \mathbf{Id}_{p\times p} \\
      \end{array}
    \right).
$$
Let $\pi$ be the following permutation of $\{1,\,\cdots,\,4p\}$ :
\begin{equation}
\pi(i)=
\left\{
\begin{array}{ll}
k+1,  &i=2k+1
  \\
p+r,  &i=2r

\\
2p+k+1, &i=2p+2k+1
\\
3p+r, &i=2r+2p
\end{array}
\right.
\end{equation}
where $k=0,\,\cdots,\,(p-1)$ and $r=1,\,\cdots,\,p$. Denote by $\mathbf{J}$ the correspond permutation matrix
$$
\mathbf{J}_{ij}=\delta_{i\,\pi(j)},\quad i,j=1,\,\cdots,\,4p.
$$
Then  $\varphi_n=\widetilde{\Psi}_n \mathbf{J}$  satisfies the following equation:
\begin{equation}\label{varphi}
\varphi_{n}=\varphi_{n-1}\left(\mathcal{A}_{\mu}\otimes A_1^{n}+\mathbf{Id}_{2p\times 2p}\otimes A_2^{n}\right),
\end{equation}
which implies that
\begin{equation*}
\varphi_{_{N}}=\Psi_{0}\mathbf{J}
\prod_{n=1}^{N}(\mathcal{A}_{\mu}\otimes A_1^n+\mathbf{Id}_{2p\times 2p}\otimes A_2^n),
\end{equation*}
and consequently that $\sigma_{N+1}^2\Psi_N^1$ satisfies equality \eqref{PsiN}.

\end{proof}
Preparing for the asymptotic study we state the following result:

\begin{Lemma}\label{bouded}
Let $(\beta_n)_{n\geq 1}$ be a  sequence of real numbers satisfying the condition \eqref{sommab_betta}. For a fixed real number $a$ let us define a sequence of $2\times 2$ matrices $(S_N(a))_{N\geq 1}$ such that:
\begin{equation}\label{S_N}
S_N(a)=\prod_{n=1}^{N-1}\left(
        \begin{array}{cc}
          a & -\beta_n \\
          -a\beta_n & 1 \\
        \end{array}
      \right)=\prod_{n=1}^{N-1}(aA_1^n+A_2^n),
\end{equation}
where $A_1^n$ and $A_1^n$ are defined by equation \eqref{A1n}.
Then

\begin{enumerate}
  \item if $|a|<1,$ $\sup\limits_{N\geq 1}\|S_N(a)\|< \infty$,
  \item if $|a|>1,$ $\sup\limits_{N\geq 1}\|(S_N(a))^{-1}\|< \infty$,
  \item if $a$ is sufficiently small, $\inf\limits_{N\geq 1}\mbox{trace}((S_N^{-1}(\displaystyle{\frac{1}{a}}))S_N(a))>0$.
\end{enumerate}

\end{Lemma}
\begin{proof}
The proof of assertions $1$ and $2$ follows directly  from the  estimates:
$$
\|aA_1^n+A_2^n\|\le 1+\beta_{n}^2\left(\frac{1+3a^2}{1- a^2}\right),\quad \quad \mbox{when} \quad |a|<1,
$$

$$
\|(aA_1^n+A_2^n)^{-1}\|\le 1+\beta_{n}^2\left(\frac{1+a^2}{a^2-1}\right),\quad \quad \mbox{when} \quad |a|>1.
$$
The proof of  assertion $3$  follows from the equality
$$ G_N(a)=\frac{1}{1-\beta_N^2}\left(
          \begin{array}{cc}
            a & -\beta_N \\
            -a\beta_N & 1 \\
          \end{array}
        \right)G_{N-1}(a)\left(
                           \begin{array}{cc}
                             a & a\beta_N \\
                             \beta_N & 1 \\
                           \end{array}
                         \right)
$$
for $G_{n}(a)=S_n^{-1}(\displaystyle{\frac{1}{a}}))S_n(a)$. Hence  $ \mbox{trace}(G_{N}(0))=\displaystyle{\frac{1}{\sigma_N^2}}$ and
condition \eqref{sommab_betta} implies that $\displaystyle{\lim_{N\rightarrow \infty}\sigma_N^2=\prod_{n=1}^{\infty}\left(1+\beta_n^2\right)<\infty}$ which achieves the proof.
\end{proof}
Actually, in the asymptotic study we work with a small value of $\mu$. Note that for a small $\mu$, matrix  $\mathcal{A}_{\mu}$ defined by \eqref{Amu} can be represented as:
$\mathcal{A}_{\mu}=\mathcal{A}_0 +\mu H$, where
\begin{equation}\label{Amusmall}
\mathcal{A}_0=\left(
              \begin{array}{cc}
                A_0^{-1} & A_0^{-1}bb^* \\
                \mathbf{0}_{p\times p} & A_0^* \\
              \end{array}
            \right) \quad H=\left(
                              \begin{array}{cc}
                               \mathbf{0}_{p\times p}  & \mathbf{0}_{p\times p} \\
                                \alpha\alpha^* & \alpha\alpha^*A_0^{-1}bb^* \\
                              \end{array}
                            \right).
\end{equation}
 Representation \eqref{Amusmall} implies that if  the spectral radius  $r(\vartheta)<1$
then there are $p$ eigenvalues  of $\mathcal{A}_{\mu}$ such that $|\lambda_i(\mu)|>1$ (in particular $\lambda_i(0),\,i=1,\,\cdots,\,p$ are the eigenvalues of $A_0^{-1}$) and
 $p$  eigenvalues of $\mathcal{A}_{\mu}$ such that $|\lambda_j(\mu)|<1,\, j=p+1,\,\cdots,\,2p.$

\begin{Lemma}\label{Laplbar}
Suppose that $r(\vartheta)<1$. Let us take $\mu=\frac{1}{N}$   and denote by
$\overline{L}_N^{\vartheta}(\mu)$:
\begin{equation}\label{lapl_eigenval}
\overline{L}_N^{\vartheta}(\mu)=\prod_{i=1}^p\left(\frac{\lambda_i(\mu)}{\lambda_i(0)}\right)^{-\frac{N}{2}}.
\end{equation}
Then, under condition \eqref{cond spec den},
\begin{equation}\label{lim 1}
\lim_{N\rightarrow \infty}\frac{L_N^{\vartheta}(\mu)}{\overline{L}_N^{\vartheta}(\mu)}=1.
\end{equation}
\end{Lemma}

\begin{proof}
Thanks  to the definition \eqref{def An} of $\mathbf{A}_n$ the equality
$$
\prod_{n=1}^{N}\det \mathbf{A}_n=\prod_{n=1}^{N}\left[(1-\beta_n^2)^p\frac{1}{\prod_{i=1}^p\lambda_i(0)}\right]
=\frac{ \left( \sigma_{N+1}^2 \right)^p}{ \prod_{i=1}^p  \lambda_i(0) ^{N} }
$$
holds. Then due to equation \eqref{Lapl_Psi} to prove  \eqref{lim 1} it is sufficient to check that
\begin{equation}\label{lim2}
\lim_{N\rightarrow \infty}\frac{\det\sigma_{N+1}^2  \Psi_N^1}{\left(\sigma_{N+1}^2\right)^p \prod_{i=1}^p[\lambda_i(\mu)]^{ N}}=1.
\end{equation}
Diagonalizing  the matrix $\mathcal{A}_{\mu}$, \textit{i.e.}, representing $\mathcal{A}_{\mu}$ as $\mathcal{A}_{\mu}=G_{\mu}D(\lambda_{i}(\mu))G_{\mu}^{-1}$ with a diagonal matrix $D(\lambda_{i}(\mu))$, we have also
$$
\mathcal{A}_{\mu}\otimes A_1^n+\mathbf{Id}_{2p\times 2p}\otimes A_2^n
$$
$$
=(G_{\mu}\otimes \mathbf{Id}_{2p\times 2p})(D(\lambda_{i}(\mu))\otimes A_1^n+\mathbf{Id}_{2p\times 2p}\otimes A_2^n)(G_{\mu}^{-1}\otimes \mathbf{Id}_{2p\times 2p}).
$$
This equation means that representation \eqref{PsiN} can be rewritten as:
\begin{equation}\label{PsiNdiag}
\sigma_{N+1}^2\Psi_N^1=\Psi_{0}\mathbf{J}(G_{\mu}\otimes \mathbf{Id}_{2p\times 2p})
 D(S_N(\lambda_{i}(\mu)))(G_{\mu}^{-1}\otimes \mathbf{Id}_{2p\times 2p}) \mathbf{J^{*}} \Psi_{0}^{*},
\end{equation}
where $D(S_N(\lambda_{i}(\mu)))$ is a block diagonal matrix with the block entries $S_N(\lambda_{i}(\mu)),\, i\le2p$ defined by equation \eqref{S_N}. 
Since $G_{0}$ is a lower triangular matrix, it follows from \eqref{PsiNdiag} that
$$
\sigma_{N+1}^2\Psi_N^1=P_{\mu}\mathcal{D}_1(S_N(\lambda_{i}(\mu)))Q_{\mu}+R_{\mu}\mathcal{D}_2(S_N(\lambda_{j}(\mu)))T_{\mu},
$$
where
$$
\lim_{\mu \rightarrow 0}P_{\mu}Q_{\mu}= \mathbf{Id}_{2p \times 2p},\quad \lim_{\mu \rightarrow 0}R_{\mu}=\mathbf{0}_{2p \times 2p},
$$
and the block diagonal matrix $\mathcal{D}_1(S_N(\lambda_{i}))$ (respectively $\mathcal{D}_2(S_N(\lambda_{j}))$) is  such that $|\lambda_i(\mu)|>1$ (respectively $|\lambda_j(\mu)|<1$ ).


Since $ \det\mathcal{D}_1(S_N(\lambda_{i}(\mu)))=\left(\sigma_{N+1}^2\right)^p\prod_{i=1}^p[\lambda_i(\mu)]^{ N}$ then, by Lemma \ref{bouded}  we get
$$
\lim_{N\rightarrow \infty}\frac{ \det\sigma_{N+1}^2 \Psi_N^1}{\det\mathcal{D}_1(S_N(\lambda_{i}(\mu)))}=1,
$$
which achieves the proof.
\end{proof}

The following statement plays a crucial role in the proofs.

\begin{Lemma}\label{lap tra}
Suppose that $r(\vartheta)<1$. Then under condition \eqref{cond spec den}, for any $\alpha \in \mathbb{R}^p$,
\begin{equation}\label{p-lap tran}
\lim_{N\rightarrow \infty}L_N^{\vartheta}(\frac{1}{N})=\exp\left(-\frac{1}{2}\alpha^*\mathcal{I}(\vartheta)\alpha\right)
\end{equation}
where  $\mathcal{I}(\vartheta)$ is the unique solution of Lyapunov equation \eqref{Lyapunov}.
\end{Lemma}

\begin{proof}
It follows immediately from Lemma~\ref{Laplbar} that under condition  \eqref{cond spec den} the limit of $ L_N^{\vartheta}(\frac{1}{N})$ does not depend on the structure of noises $\xi$, \textit{ i.e.}, does not depend on $\beta_{n}$. Thus, this limit is the same as for the classical i.i.d. situation, when $\beta_n=0,n \ge 1$.

\end{proof}

\begin{rem}\label{Eqlayp_eignv}
It is worth mentioning that  equation \eqref{lapl_eigenval} says that
$$
\sum_{i=1}^{p}\frac{d}{d\mu}\ln(\lambda_{i}(\mu))|_{\mu=0}= \alpha^*\mathcal{I}(\vartheta)\alpha,
$$
where $\lambda_i(\mu)$ are the eigenvalues of $\mathcal{A}_{\mu}$ such that $|\lambda_i(\mu)|>1$  and $\mathcal{I}(\vartheta)$ is the solution of the Lyapunov equation \eqref{Lyapunov}. Of course, this equality  can be proved independently.
\end{rem}

\section{Proofs}\label{proofs}
\subsection{Proof of Theorem ~\ref{p-dimension}}
The statement of Theorem follows from Lemma~\ref{lap tra} since \eqref{p-lap tran} implies immediately that
\begin{equation}\label{con in law}
\mathbf{P}_{\vartheta}-\lim_{N\rightarrow \infty}\frac{1}{N}\langle M\rangle_N=\mathcal{I}(\vartheta),
\end{equation}
and, hence also due to the central limit theorem for martingales,
$$
\frac{1}{\sqrt{N}}M_N \overset{\textit{law}}{\Rightarrow} \mathcal{N}\left(\mathbf{0},\,\mathcal{I}(\vartheta)\right).
$$

\subsection{Proof of Theorem ~\ref{strong}}
Due to the strong law of large numbers for martingales, in order to proof the strong consistency we have only to check  that
$$
\lim_{N \rightarrow \infty} \langle M\rangle_N =+\infty    \quad a.s.,
$$
or, equivalently that for a one  fixed constant $\mu>0$
\begin{equation}\label{eq strong}
\lim_{N\rightarrow \infty}\mathbf{E}_{\vartheta}\exp\left(-\frac{\mu}{2}\langle M\rangle_N\right)=0.
\end{equation}
But in the case when $p=1$ the ingredients in the right hand side of formulas \eqref{Lapl_Psi}- \eqref{PsiN} with $\alpha=1$ can be given more explicitly:

$$
\prod_{n=1}^{N}\det \mathbf{A}_n=\vartheta^{N}\sigma_{N+1}^2,
$$
and
\begin{equation}\label{psi_p1}
\sigma_{N+1}^2\Psi_N^1=\frac{1-\lambda_{-}}{\lambda_{+} - \lambda_{-}}S_{N}(\frac{\lambda_{+}}{\vartheta})+ \frac{\lambda_{+}-1}{\lambda_{+} - \lambda_{-}}S_{N}(\frac{\lambda_{-}}{\vartheta}),
\end{equation}
where the matrix $S_{N}(a)$ is defined by equation \eqref{S_N},
$$
\frac{\lambda_{\pm}}{\vartheta}=\frac{\vartheta^2+\mu+1\pm \sqrt{(\mu+(1-\vartheta)^2)(\mu+(1+\vartheta)^2)}}{2\vartheta}
$$
are the two eigenvalues of the matrix $
\displaystyle{\mathcal{A}_{\mu}=\left(
                                         \begin{array}{cc}
                                           \frac{1}{\vartheta} & \frac{1}{\vartheta} \\
                                           \mu\frac{1}{\vartheta} & \mu\frac{1}{\vartheta}+\vartheta \\
                                         \end{array}
                                       \right)}
$
Note  that $\displaystyle{\frac{\lambda_{+}}{\vartheta}\frac{\lambda_{-}}{\vartheta}=1,\, \left|\frac{\lambda_{+}}{\vartheta}\right|>1}$ and  $\lambda_{+}>1$ for every $\mu>0$ and $\vartheta \in \mathbb{R}$.
Equations \eqref{psi_p1} and \eqref{S_N} imply that for $\kappa=\frac{\lambda_{+}-1}{1-\lambda_{-}}$
$$
\det \Psi_N^1=\left(\frac{1}{\sigma^2_{N+1}}\right)^2 \left(\frac{1-\lambda_{-}}{\lambda_{+} - \lambda_{-}}\right)^{2}\det\left(S_{N}(\frac{\lambda_{+}}{\vartheta})\right)
\det\left(\mathbf{Id}_{2\times 2}+ \kappa(S_{N}(\frac{\lambda_{+}}{\vartheta}))^{-1}S_{N}(\frac{\lambda_{-}}{\vartheta})\right)
$$
and that
 $$
 \det\left(S_{N}(\frac{\lambda_{+}}{\vartheta})\right)=\vartheta^{-N}\lambda_{+}^{N}\sigma_{N+1}^2.
 $$
Thanks to Lemma  \ref{bouded} 

$$
\det\left(\mathbf{Id}_{2\times 2}+ \kappa(S_{N}(\frac{\lambda_{+}}{\vartheta}))^{-1}S_{N}(\frac{\lambda_{-}}{\vartheta})\right)
$$
is uniformly bounded and separated from $0$ when $\mu$ is sufficiently large (and so $a=\frac{\lambda_{-}}{\vartheta}$ is sufficiently small). Since  $\lambda_{+}>1$, we obtain that
$$
\lim_{N\rightarrow \infty}L_N^{\vartheta}(\mu)=c\lim_{N\rightarrow \infty}\lambda_{+}^{-\frac{N}{2}}=0.
$$
The uniform
consistency  and the uniform convergence of the moments on compacts  $\m{K} \subset (-1,\,1)$ follow from the estimates
(see \cite{Liptser2}, Eq.17.51):
$$
\mathbf{E}_{\vartheta}\left(\frac{1}{N}\langle M\rangle_N\right)^{-q} \leq (1-\vartheta^2)^{-q},
$$

$$
\mathbf{E}\left(\frac{1}{\sqrt{N}}M_N\right)^q \leq \left(\sqrt{1-\vartheta^2}\right)^q.
$$

\begin{rem}
It is worth mentioning that even in a stationary autoregressive  models of order $1$ with strongly dependent noises
the  Least Square Estimator
$ \widetilde{\vartheta}_N = \frac{\sum_{n=1}^N X_{n-1} X_n }{ \sum_{n=1}^N X^2_{n-1}}$
is not  consistent.
\end{rem}

\bibliographystyle{plain}
\bibliography{recherchebib}

\end{document}